# INVARIANT TRANSPORTS OF STATIONARY RANDOM MEASURES AND MASS-STATIONARITY

By Günter Last and Hermann Thorisson

*Universität Karlsruhe and University of Iceland*

We introduce and study invariant (weighted) transport-kernels balancing stationary random measures on a locally compact Abelian group. The first main result is an associated fundamental invariance property of Palm measures, derived from a generalization of Neveu's exchange formula. The second main result is a simple sufficient and necessary criterion for the existence of balancing invariant transport-kernels. We then introduce (in a nonstationary setting) the concept of mass-stationarity with respect to a random measure, formalizing the intuitive idea that the origin is a typical location in the mass. The third main result of the paper is that a measure is a Palm measure if and only if it is mass-stationary.

**1. Introduction.** We consider (jointly) stationary random measures on a locally compact Abelian group $G$, for instance, $G = \mathbb{R}^d$. A *transport-kernel* is a Markovian kernel $T$ that distributes mass over $G$ and depends on both $\omega$ in the underlying sample space $\Omega$ and a location $s \in G$. The number $T(\omega, s, B)$ is the proportion of mass *transported* from location $s$ to the set $B$. More generally, a *weighted* transport-kernel is a kernel $T$ which need not be Markovian. If $T$ is finite, then the mass at $s$ is weighted by $T(\omega, s, G)$ before being transported by the normalized $T$. In general, we assume that $T(\omega, s, B)$ is finite for compact $B$ but allow that $T(\omega, s, G) = \infty$. A kernel $T$ is *invariant* if it is invariant under joint shifts of all three arguments. If $\xi$ and $\eta$ are random measures on $G$ such that $\xi T = \eta$, then $T$ is $(\xi, \eta)$-*balancing* and, in particular, if $\xi = \eta$, then $T$ is $\xi$-*preserving*. Sometimes an invariant $T$ can be reduced to an *allocation rule* $\tau$ (depending on $\omega \in \Omega$) that maps each location $s$ to a new location $\tau(s)$ in an invariant way. In fact, we might think of an invariant transport-kernel $T$ as the conditional distribution of a randomized allocation rule.

The aim of this paper is to treat three interwoven aspects of invariant weighted transport-kernels: basic invariance properties of Palm measures are presented in Sections 3 and 4, a general existence result in Section 5, and an intrinsic characterization of Palm measures—which we call *mass-stationarity*—in Sections 6 and 7. Below we sketch these results against their background.

INVARIANCE PROPERTIES (SECTIONS 3 AND 4). It is a fundamental and classical theorem that the Palm distribution of a stationary point process on the line is invariant under shifts to the next point on the right. In fact, there is a unique correspondence between such stationary point processes and stationary sequences of interpoint distances; see Theorem 11.4 in [12] and the references given there.

In higher dimensions, the situation is more complicated. Mecke [21] found an intrinsic characterization of Palm measures using an integral equation; see (2.7) below. Geman and Horowitz [5] showed that a stationary random measure is distributionally invariant under shifts associated with allocation rules preserving Haar measure. Independently, Mecke [22] derived more general invariance properties of Palm measures; see Remark 4.6. In the point process case one of his results was rediscovered in [27]; see also [7]. Port and Stone [24] studied what they called (translation-invariant) *marked motion process*. In particular, they derived a certain transport property of Palm distributions (called *tagged particle distributions*) associated with such processes; see Example 4.3 below. Holroyd and Peres [11] considered (in case $G = \mathbb{R}^d$) invariant transport-kernels $T$ balancing Lebesgue measure and an ergodic point process $\eta$ of intensity 1. [They call the random measure $T(s, dt)\, ds$ on $\mathbb{R}^d \times \mathbb{R}^d$ a *transport rule*.] Theorem 16 in that paper states that if the origin 0 is shifted to a randomized location with conditional distribution $T(\omega, 0, \cdot)$, then the stationary distribution of $\eta$ is transformed into the Palm distribution of $\eta$; this is an example of *shift-coupling*. An analogous result holds for discrete groups.

Neveu's [23] well-known exchange formula (see Remark 3.7) is an apparently quite different property of Palm measures. We will generalize Neveu's result in our Theorem 3.6. This is then actually the key to obtaining the general invariance property of Theorem 4.1, containing all the invariance results mentioned above. Another crucial idea for Theorem 4.1 is that any balancing invariant weighted transport-kernel has an *inverse* invariant transport-kernel.

EXISTENCE (SECTION 5). Liggett [19] constructed an allocation rule, transporting counting measures on the integers to the Bernoulli (1/2) random measure with intensity 1. (He also treated a general Bernoulli parameter $p$ and the Poisson process on the line.) Triggered by Liggett's paper, invariant transports of the Lebesgue measure on $\mathbb{R}^d$ (or of the counting measure on $\mathbb{Z}^d$) to an ergodic point process of intensity 1 have received considerable attention in recent years (see [3, 9, 11, 14]). In particular, Holroyd and Peres [11] constructed an explicit algorithm based on a so-called stable marriage allocation.

Actually, an abstract group-coupling result from [26] implies that shift-couplings exist in the above cases; see [27]. In Section 5 we shall apply that result to stationary random measures with finite intensities to prove that there exists an invariant balancing transport-kernel if and only if the two random measures have the same intensity conditional on the invariant $\sigma$-field.

MASS-STATIONARITY (SECTIONS 6 AND 7). Thorisson [27] calls a simple point process on $\mathbb{R}^d$ *point-stationary* if (loosely speaking) it looks distributionally the same from all its points, just like stationarity means that the process looks distributionally the same from all locations in $\mathbb{R}^d$. The formal definition in that paper required a joint distributional invariance under shifts associated with certain preserving *randomized* allocation rules. The main result was that point-stationarity is a characterizing property of Palm versions of stationary point processes. The question of whether the external randomization could be removed from the definition inspired considerable research activity; see [4, 10, 14, 29]. Finally, Heveling and Last [7, 8] showed that this can be done.

In Section 6 we extend the concept of point-stationarity to a random measure $\xi$. We will call $\xi$ *mass-stationary* if (again loosely speaking) it looks distributionally the same from all locations in its mass. The formal definition of this property is quite subtle and requires some joint distributional invariance; see Definition 6.1 and Remark 6.2. Our Theorem 6.3 states that mass-stationarity is a characterizing property of Palm versions of stationary random measures. In Section 7 we discuss mass-stationarity briefly. In particular, we show that mass-stationarity is equivalent to distributional invariance under bounded invariant $\xi$-preserving weighted transport-kernels. On the other hand, Example 7.1 shows that invariance under preserving (nonrandomized) allocation rules is not enough to imply mass-stationarity. We conclude with five open problems.

**2. Preliminaries on stationary random measures.** We choose to work in the abstract setting of a flow acting on the underlying sample space (see [5, 22, 23]), and with $\sigma$-finite measures rather than with probability measures; see Remark 2.6.

We consider a topologial Abelian group $G$ that is assumed to be a locally compact, second countable Hausdorff space with Borel $\sigma$-field $\mathcal{G}$. On $G$ there exists an invariant measure $\lambda$, that is unique up to normalization. A measure $\mu$ on $G$ is *locally finite* if it is finite on compact sets. We denote by $\mathbf{M}$ the set of all locally finite measures on $G$, and by $\mathcal{M}$ the cylindrical $\sigma$-field on $\mathbf{M}$ which is generated by the evaluation functionals $\mu \mapsto \mu(B)$, $B \in \mathcal{G}$. The *support* $\operatorname{supp}\mu$ of a measure $\mu \in \mathbf{M}$ is the smallest closed set $F \subset G$ such that $\mu(G \setminus F) = 0$. By $\mathbf{N} \subset \mathbf{M}$, we denote the measurable set of all (simple) *counting measures* on $G$, that is, the set of all those $\mu \in \mathbf{M}$ with discrete support and $\mu\{s\} := \mu(\{s\}) \in \{0,1\}$ for all $s \in G$. We can and will identify $\mathbf{N}$ with the class of all *locally finite* subsets of $G$, where a set is called locally finite if its intersection with any compact set is finite.

In this paper we mostly work on a $\sigma$-finite measure space $(\Omega, \mathcal{F}, \mathbb{P})$ (but see also Remark 2.6). However, we will consider several measures on $(\Omega, \mathcal{F})$. A *random measure* on $G$ is a measurable mapping $\xi : \Omega \to \mathbf{M}$ and a (simple) *point process* on $G$ is a measurable mapping $\xi : \Omega \to \mathbf{N}$. A random measure $\xi$ can also be regarded as a *kernel* from $\Omega$ to $G$. Accordingly, we write $\xi(\omega, B)$ instead of $\xi(\omega)(B)$. If $\xi$ is a random measure, then the mapping $(\omega, s) \mapsto \mathbf{1}\{s \in \operatorname{supp}\xi(\omega)\}$ is measurable.

We assume that $(\Omega, \mathcal{F})$ is equipped with a *measurable flow* $\theta_s : \Omega \to \Omega$, $s \in G$. This is a family of measurable mappings such that $(\omega, s) \mapsto \theta_s \omega$ is measurable, $\theta_0$ is the identity on $\Omega$ and

$$\theta_s \circ \theta_t = \theta_{s+t}, \qquad s, t \in G, \tag{2.1}$$

where 0 denotes the neutral element in $G$ and $\circ$ denotes composition. A random measure $\xi$ on $G$ is called *invariant* (or *flow-adapted*) if

$$\xi(\theta_s \omega, B - s) = \xi(\omega, B), \qquad \omega \in \Omega, \ s \in G, \ B \in \mathcal{G}. \tag{2.2}$$

A measure $\mathbb{P}$ on $(\Omega, \mathcal{F})$ is called *stationary* if it is invariant under the flow, that is,

$$\mathbb{P} \circ \theta_s = \mathbb{P}, \qquad s \in G,$$

where $\theta_s$ is interpreted as a mapping from $\mathcal{F}$ to $\mathcal{F}$ in the usual way:

$$\theta_s A := \{\theta_s \omega : \omega \in A\}, \qquad A \in \mathcal{F}, \ s \in G.$$

Because of the next examples we may think of $\theta_s \omega$ as of $\omega$ *shifted* by $-s$.

EXAMPLE 2.1. Consider the measurable space $(\mathbf{M}, \mathcal{M})$ and define for $\mu \in \mathbf{M}$ and $s \in G$ the measure $\theta_s \mu$ by $\theta_s \mu(B) := \mu(B + s)$, $B \in \mathcal{G}$. Then $\{\theta_s : s \in G\}$ is a measurable flow and the identity $\xi$ on $\mathbf{M}$ is an invariant random measure. A stationary probability measure on $(\mathbf{M}, \mathcal{M})$ can be interpreted as the distribution of a *stationary random measure*. Let $(E, \mathcal{E})$ be some measurable space and denote by $\mathbf{M}_E$ the set of all measures $\mu$ on $G \times E$ such that $\mu(\cdot \times E) \in \mathbf{M}$. Let $\mathcal{M}_E$ be the $\sigma$-field on $\mathbf{M}_E$ generated by the mappings $\mu \mapsto \mu(B)$, $B \in \mathcal{G} \otimes \mathcal{E}$. For $\mu \in \mathbf{M}_E$ and $s \in G$, let $\theta_s \mu$ be the measure $\mu \in \mathbf{M}_E$ satisfying $\theta_s \mu(B \times C) := \mu((B+s) \times C)$, for all $B \in \mathcal{G}$ and $C \in \mathcal{E}$. Then $\{\theta_s : s \in G\}$ is a measurable flow on $\mathbf{M}_E$. A stationary probability measure on $(\mathbf{M}_E, \mathcal{M}_E)$ can be interpreted as the distribution of a *stationary marked random measure*.

REMARK 2.2. Since a random measure $\xi$ is a random element in $\mathbf{M}$, we can rewrite the invariance condition (2.2) as

$$\xi(\theta_s \omega) = \theta_s \xi(\omega), \qquad \omega \in \Omega,$$

where we use $\theta_s$, $s \in G$, to denote both the abstract flow and the specific flow defined in Example 2.1. Therefore, (2.2) is also referred to as *flow-covariance*. We follow here the terminology of [13]. A similar remark applies to invariant weighted transport-kernels, to be defined below.

EXAMPLE 2.3. Let $(E, \mathcal{E})$ be a Polish space and assume that $\Omega$ is the space of all measures $\omega$ on $G \times E \times G \times E$ such that $\omega(B \times E \times G \times E)$ and $\omega(G \times E \times B \times E)$ are finite for compact $B \subset G$. The $\sigma$-field $\mathcal{F}$ is defined analogously as in Example 2.1. It is stated in [24] (and can be proved as in [20]) that $(\Omega, \mathcal{F})$ is a Polish space. For $s \in G$ and $\omega \in \Omega$, we let $\theta_s \omega$ denote the measure satisfying

$$\theta_s \omega(B \times C \times B' \times C') = \omega((B+s) \times C \times (B'+s) \times C')$$

for all $B, B' \in \mathcal{G}$ and $C, C' \in \mathcal{E}$. The random measures $\xi$ and $\eta$ defined by $\xi(\omega, \cdot) := \omega(\cdot \times E \times G \times E)$ and $\eta(\omega, \cdot) := \omega(G \times E \times \cdot \times E)$ are invariant. Port and Stone [24] (see also [6]) call a stationary probability measure on $(\Omega, \mathcal{F})$ concentrated on the set of integer-valued $\omega \in \Omega$ a (translation invariant) *marked motion process*. The idea is that the (marked) points of $\xi$ move to the points of $\eta$ in one unit of time.

EXAMPLE 2.4. Assume that $(\Omega, \mathcal{F}) = (G, \mathcal{G})$ and $\theta_s \omega := \omega + s$. Then the Haar measure $\mathbb{P} := \lambda$ is stationary. If $\mathbb{P}'$ is a probability measure on $G$, then $\xi(\omega, B) := \mathbb{P}'(B + \omega)$ defines an invariant random measure.

Let $\mathbb{P}$ be a stationary $\sigma$-finite measure on $(\Omega, \mathcal{F})$ and $\xi$ an invariant random measure on $G$. Then $\xi$ is stationary in the usual sense, that is, $\mathbb{P}(\xi \in \cdot) = \mathbb{P}(\theta_s \xi \in \cdot)$ for all $s \in G$, where we have used the notation of Example 2.1. Let $B \in \mathcal{G}$ be a set with positive and finite Haar measure $\lambda(B)$. The measure

$$(2.3) \quad \mathbb{P}_\xi(A) := \lambda(B)^{-1} \iint \mathbf{1}_A(\theta_s \omega) \mathbf{1}_B(s) \xi(\omega, ds) \mathbb{P}(d\omega), \qquad A \in \mathcal{F},$$

is called the *Palm measure* of $\xi$ (with respect to $\mathbb{P}$); see [21]. This measure is $\sigma$-finite. As the definition (2.3) is independent of $B$, we can use a monotone class argument to conclude the *refined Campbell theorem*

$$\iint f(\theta_s\omega,s)\xi(\omega,ds)\mathbb{P}(d\omega) = \iint f(\omega,s)\,ds\,\mathbb{P}_\xi(d\omega)$$

for all measurable $f:\Omega \times G \to [0,\infty)$, where $ds$ refers to integration with respect to the Haar measure $\lambda$. Using a standard convention in probability theory, we write this as

$$(2.4) \qquad \mathbb{E}_\mathbb{P}\bigg[\int f(\theta_s,s)\xi(ds)\bigg] = \mathbb{E}_{\mathbb{P}_\xi}\bigg[\int f(\theta_0,s)\,ds\bigg],$$

where $\mathbb{E}_\mathbb{P}$ and $\mathbb{E}_{\mathbb{P}_\xi}$ denote integration with respect to $\mathbb{P}$ and $\mathbb{P}_\xi$, respectively.

EXAMPLE 2.5. Consider the setting of Example 2.3. Let $\tilde{\xi}$ be the random element in $\mathbf{M}_E$ (cf. Example 2.1) defined by $\tilde{\xi}(\omega,\cdot) := \omega(\cdot \times G \times E)$. Assume that $\mathbb{P}(\tilde{\xi} \in \cdot)$ is $\sigma$-finite. Then there is a Markov kernel $K$ from $\mathbf{M}_E$ to $\Omega$ satisfying

$$\mathbb{P} = \int K(\mu,\cdot)\mathbb{P}(\tilde{\xi} \in d\mu).$$

By Theorem 3.5 in [13], we can assume that $K$ is invariant in the sense that

$$K(\theta_s\mu,\theta_sA) = K(\mu,A), \qquad s \in G,\ \mu \in \mathbf{M}_E,\ A \in \mathcal{F},$$

where $\{\theta_s\}$ denotes the flow on both $\mathbf{M}_E$ and $\Omega$. Of course, if $\mathbb{P}$ is a probability measure, then $K(\tilde{\xi},A)$ is a version of the conditional probability of $A \in \mathcal{F}$ given $\tilde{\xi}$. Using invariance of $K$, it is straightforward to check that

$$(2.5) \qquad \mathbb{P}_\xi = \int K(\mu,\cdot)\mathbb{P}_\xi(\tilde{\xi} \in d\mu).$$

If the *intensity* $\mathbb{P}_\xi(\Omega)$ of $\xi$ is positive and finite, then the normalized Palm measure

$$\mathbb{P}_\xi^0 := \mathbb{P}_\xi(\Omega)^{-1}\mathbb{P}_\xi$$

is called the *Palm probability measure* of $\xi$ (w.r.t. $\mathbb{P}$). Note that $\mathbb{P}_\xi$ and $\mathbb{P}_\xi^0$ are both defined on the underlying space $(\Omega,\mathcal{F})$. The stationary measure $\mathbb{P}$ can be recovered from the Palm measure $\mathbb{P}_\xi$ using a measurable function $\tilde{h}:\mathbf{M} \times G \to [0,\infty)$ satisfying $\int \tilde{h}(\mu,s)\mu(ds) = 1$, whenever $\mu \in \mathbf{M}$ is not the null measure. For one example of such a function we refer to [21]. We then have the *inversion formula*

$$(2.6) \qquad \mathbb{E}_\mathbb{P}[\mathbf{1}\{\xi(G) > 0\}f] = \mathbb{E}_{\mathbb{P}_\xi}\bigg[\int \tilde{h}(\xi \circ \theta_{-s},s)f(\theta_{-s})\,ds\bigg],$$

for all measurable $f:\Omega \to [0,\infty)$; see Satz 2.4 in Mecke [21]. This is a direct consequence of the refined Campbell theorem (2.4).

Let $\xi$ be an invariant random measure and $\mathbb{Q}$ be a $\sigma$-finite measure on $(\Omega,\mathcal{F})$ satisfying $\mathbb{Q}(\xi(G)=0)=0$. Satz 2.5 in Mecke [21] says that there is $\sigma$-finite stationary measure $\mathbb{P}$ on $(\Omega,\mathcal{F})$ such that $\mathbb{Q}$ is a Palm measure $\mathbb{P}_\xi$ of $\xi$ with respect to $\mathbb{P}$ if and only if for all measurable $g:\Omega \times G \to [0,\infty)$

$$(2.7) \qquad \mathbb{E}_\mathbb{Q}\left[\int g(\theta_s, -s)\xi(ds)\right] = \mathbb{E}_\mathbb{Q}\left[\int g(\theta_0, s)\xi(ds)\right].$$

Mecke proved his result in the canonical framework of Example 2.1. But his proof applies in our more general framework as well. The necessity of (2.7) is a special case of Neveu's exchange formula; see Remark 3.7. To prove that (2.7) is also sufficient for $\mathbb{Q}$ to be a Palm measure, one can use the function $\tilde{h}$ in (2.6) to define a $\sigma$-finite measure $\mathbb{P}$ by

$$\mathbb{P}(A) := \mathbb{E}_\mathbb{Q}\left[\int \tilde{h}(\xi \circ \theta_{-s}, s)\mathbf{1}_A(\theta_{-s})\,ds\right], \qquad A \in \mathcal{F}.$$

It can be shown, as in [21], that (2.7) implies stationarity of $\mathbb{P}$ and $\mathbb{Q} = \mathbb{P}_\xi$.

REMARK 2.6. We would like to mention two reasons (other than just generality) why we are not assuming the stationary measure $\mathbb{P}$ to be a probability measure. First, some of the fundamental results can be more easily stated this way. An example is the one-to-one correspondence between $\mathbb{P}$ and the Palm measure $\mathbb{P}_\xi$ (see [21]). Otherwise, extra technical integrability assumptions are required (see, e.g., Theorem 11.4 in [12]). A second reason is that in some applications it is the Palm probability measure that has a probabilistic interpretation (see, e.g., [30]). This measure can be defined whenever the (stationary) intensity is positive and finite; see Example 3.5 for a simple illustration of this fact.

**3. Transport-kernels and an exchange formula.** A *transport-kernel* (on $G$) is a Markovian kernel $T$ from $\Omega \times G$ to $G$. It is helpful to think of $T(\omega, s, B)$ as the proportion of mass transported from location $s$ to the set $B$, when $\omega$ is given. A *weighted transport-kernel* is a kernel $T$ from $\Omega \times G$ to $G$ such that $T(\omega, s, \cdot)$ is locally finite for all $(\omega, s) \in \Omega \times G$. A weighted transport-kernel $T$ is called *invariant* if

$$(3.1) \quad T(\theta_t\omega, s-t, B-t) = T(\omega, s, B), \qquad s, t \in G,\ \omega \in \Omega,\ B \in \mathcal{G}.$$

This is equivalent to $T(\theta_t\omega, 0, B-t) = T(\omega, t, B)$ for all $t$, $\omega$ and $B$. Quite often we use the short-hand notation $T(s, \cdot) := T(\theta_0, s, \cdot)$. If $\tilde{T}$ is kernel from $\Omega$ to $G$ such that $\tilde{T}(\omega, \cdot)$ is locally finite for all $\omega \in \Omega$, then $T(\omega, s, B) := \tilde{T}(\theta_s\omega, B-s)$ defines an invariant weighted transport-kernel $T$ on $G$.

REMARK 3.1. Let $T$ be a weighted transport-kernel on $G$ and $\xi$ an invariant random measure on $G$. Assume that $\eta := \int T(\omega, s, \cdot)\xi(\omega, ds)$ is locally finite for each $\omega \in \Omega$. If $\xi$ and $T$ are invariant, then $\eta$ is invariant too. More generally, the random measure $\psi$ on $G \times G$, defined by $\psi(d(s,t)) := T(s, dt)\xi(ds)$, is invariant in the obvious way. Assume that $\mathbb{P}$ is a stationary probability measure on $(\Omega, \mathcal{F})$. Then $\psi$ is generalizing the *marked motion processes* of [24]; see Example 2.3. Another special case of $\psi$ are the *transport rules* of [11] arising in case $\xi$ is the Lebesgue measure on $G = \mathbb{R}^d$. While our terminology was motivated by [24] and [11] (and Example 3.2 below), we found it more convenient to put the focus on the kernel $T$. Our interpretation is that $T$ *transports* $\xi$ to $\eta$ in an invariant way.

EXAMPLE 3.2. Consider a measurable function $\kappa: \Omega \times G \times G \to [0, \infty)$ and assume that $\kappa$ is invariant, that is,

$$(3.2) \qquad \kappa(\theta_r \omega, s - r, t - r) = \kappa(\omega, s, t), \qquad \omega \in \Omega, \ r, s, t \in G.$$

Let $\eta$ be an invariant random measure on $G$ and define

$$T(s, B) := \int_B \kappa(s, t)\eta(dt).$$

Then (3.1) holds. Such functions $\kappa$ occur in the *mass-transport principle*; see [2] and Remark 3.8 below. If $\eta$ is a simple point process and $t \in \operatorname{supp} \eta(\omega)$, then the number $\kappa(\omega, s, t)$ is interpreted as the mass sent from $s$ to $t$ when the configuration $\omega$ is given.

Let $\xi$ and $\eta$ be two invariant random measures on $G$. A weighted transport-kernel $T$ on $G$ is called $(\xi, \eta)$-*balancing*, if $T$ transports $\xi$ to $\eta$, that is, if for all $\omega \in \Omega$

$$(3.3) \qquad \int T(\omega, s, \cdot)\xi(\omega, ds) = \eta(\omega, \cdot).$$

In case $\xi = \eta$ we also say that $T$ is $\xi$-*preserving*. If $\mathbb{Q}$ is a measure on $(\Omega, \mathcal{F})$ such that (3.3) holds for $\mathbb{Q}$-a.e. $\omega \in \Omega$, then we say that $T$ is $\mathbb{Q}$-a.e. $(\xi, \eta)$-balancing.

EXAMPLE 3.3. Consider the setting of Example 2.3 and let $\mathbb{P}$ be a $\sigma$-finite stationary measure on $(\Omega, \mathcal{F})$. Then there is a $\mathbb{P}$-a.e. $(\xi, \eta)$-*balancing* invariant transport-kernel $T$. To see this, we define a measure $M$ on $\Omega \times G \times G$ by

$$M := \iiint \mathbf{1}\{(\omega, s, t) \in \cdot\}\omega(ds \times E \times dt \times E)\mathbb{P}(d\omega).$$

Stationarity of $\mathbb{P}$ implies that

$$(3.4) \qquad \int \mathbf{1}\{(\theta_r \omega, s - r, t - r) \in \cdot\}M(d(\omega, s, t)) = M, \qquad r \in G.$$

The measure

$$M' := M(\cdot \times G) = \iint \mathbf{1}\{(\omega, s) \in \cdot\} \xi(\omega, ds) \mathbb{P}(d\omega)$$

is $\sigma$-finite. Hence, we can apply Theorem 3.5 in Kallenberg [13] to obtain an invariant transport-kernel $T$ satisfying $M(d(\omega, s, t)) = T(\omega, s, dt) M'(d(\omega, s))$. [In fact, the theorem yields an invariant kernel $T'$, satisfying this equation. But in our specific situation we have $T'(\omega, s, G) = 1$ for $M'$-a.e. $(\omega, t)$, so that $T'$ can be modified in an obvious way to yield the desired $T$.] It is easy to see that $T$ is indeed $\mathbb{P}$-a.e. $(\xi, \eta)$-balancing.

If $\mathbb{P}$ is concentrated on the set $\Omega'$ of all integer-valued $\omega \in \Omega$, then one possible choice of a $\mathbb{P}$-a.e. $(\xi, \eta)$-balancing invariant transport-kernel $T$ (consistent with the above proof) is

$$(3.5) \quad T(\omega, s, \cdot) = \frac{1}{\xi(\omega, \{s\})} \sum_{t: \omega(s,t) > 0} \omega(s, t) \delta_t, \qquad \omega \in \Omega',$$

if $\xi(\omega, \{s\}) > 0$ and $T(\omega, s, \cdot) := \delta_s$ otherwise, where $\omega(s, t) := \omega(\{s\} \times E \times \{t\} \times E)$. A general criterion for the existence of balancing transport-kernels is given in Section 5.

EXAMPLE 3.4. Consider a measurable mapping $\tau: \Omega \times G \to G$. Then the transport-kernel $T(s, \cdot) := \delta_{\tau(s)}$ is invariant if and only if $\tau$ is *covariant* in the sense that

$$\tau(\theta_t \omega, s - t) = \tau(\omega, s) - t, \qquad s, t \in G, \ \omega \in \Omega.$$

In this case, following [11], we call $\tau$ an *allocation rule*. Covariance of $\tau$ is equivalent to $\tau(\theta_t \omega, 0) = \tau(\omega, t) - t$ for all $\omega, t$. Writing $\pi := \tau(0) := \tau(\theta_0, 0)$, we can express this as $\tau(s) = \pi \circ \theta_s + s$. Any measurable mapping $\pi: \Omega \to G$ can be used to generate an allocation rule this way. We interpret an allocation rule as allocating (or transporting), given $\omega \in \Omega$, an actual unit of mass close to $s$ to a new location $\tau(\omega, s)$. Let $\xi$ and $\eta$ be two random measures on $G$. The transport-kernel $T(s, \cdot) := \delta_{\tau(s)}$ is $(\xi, \eta)$-balancing iff

$$(3.6) \quad \int \mathbf{1}\{\tau(s) \in \cdot\} \xi(ds) = \eta.$$

We then say that $\tau$ is $(\xi, \eta)$-*balancing*. In case $\xi = \eta$ we also say that $\tau$ is $\xi$-*preserving*.

EXAMPLE 3.5. Consider the setting of Example 2.4. Letting $\mathbb{P}'$ and $\xi$ be as in that example, we obtain from an easy calculation that $\mathbb{P}_\xi = \mathbb{P}'$.

Now let $K$ be a Markovian kernel from $G$ to $G$ and define $\mathbb{P}'' := \int K(s, \cdot) \mathbb{P}'(ds)$. The probability measure $\iint \mathbf{1}\{(s, t) \in \cdot\} K(s, dt) \mathbb{P}'(ds)$ is a *coupling* of $\mathbb{P}'$ and $\mathbb{P}''$. In the Monge–Kantorovich mass transportation theory (see, e.g., [25]) $K$ is interpreted as transporting the mass distribution $\mathbb{P}'$ to $\mathbb{P}''$. Let $\eta$ be the invariant random measure $\eta(\omega, B) := \mathbb{P}''(B + \omega)$, and define an invariant transport-kernel $T$ by $T(\omega, s, B) := K(\omega + s, B + \omega)$. Then

$$\int T(\omega, s, B)\xi(\omega, ds) = \int K(s, B + \omega)\mathbb{P}'(ds) = \mathbb{P}''(B + \omega) = \eta(\omega, B),$$

that is, $T$ is $(\xi, \eta)$-balancing. Conversely, if $\mathbb{P}'$ and $\mathbb{P}''$ are given, and $T$ is a $(\xi, \eta)$-balancing transport-kernel (with $\xi$ and $\eta$ defined as before), then $T(0, s, B)$ is a Markovian kernel transporting $\mathbb{P}'$ to $\mathbb{P}''$.

We now prove an important transport property of Palm measures.

THEOREM 3.6. *Let $\mathbb{P}$ be a $\sigma$-finite stationary measure on $(\Omega, \mathcal{F})$. Consider two invariant random measures $\xi$ and $\eta$ on $G$ and let $T$ and $T^*$ be invariant weighted transport-kernels satisfying for $\mathbb{P}$-a.e. $\omega \in \Omega$*

$$\iint \mathbf{1}\{(s,t) \in \cdot\}T(\omega, s, dt)\xi(\omega, ds) \tag{3.7}$$
$$= \iint \mathbf{1}\{(s,t) \in \cdot\}T^*(\omega, t, ds)\eta(\omega, dt).$$

*Then we have for any measurable function $h: \Omega \times G \to [0, \infty)$ that*

$$\mathbb{E}_{\mathbb{P}_\xi}\left[\int h(\theta_t, -t)T(0, dt)\right] = \mathbb{E}_{\mathbb{P}_\eta}\left[\int h(\theta_0, t)T^*(0, dt)\right]. \tag{3.8}$$

PROOF. Let $B \in \mathcal{G}$ satisfy $\lambda(B) = 1$ and take a measurable $h: \Omega \times \mathbb{R}^d \to [0, \infty)$. From the definition (2.3) of $\mathbb{P}_\xi$ and (2.1) we obtain

$$I := \mathbb{E}_{\mathbb{P}_\xi}\left[\int h(\theta_t, -t)T(\theta_0, 0, dt)\right]$$
$$= \mathbb{E}_{\mathbb{P}}\left[\iint \mathbf{1}_B(s)h(\theta_{s+t}, -t)T(\theta_s, 0, dt)\xi(ds)\right]$$
$$= \mathbb{E}_{\mathbb{P}}\left[\iint \mathbf{1}_B(s)h(\theta_t, -t + s)T(\theta_0, s, dt)\xi(ds)\right],$$

where we have used the invariance (3.1) to get the second equation. Now we can apply assumption (3.7) to get

$$I = \mathbb{E}_{\mathbb{P}}\left[\iint \mathbf{1}_B(s)h(\theta_t, s - t)T^*(\theta_0, t, ds)\eta(dt)\right]$$
$$= \mathbb{E}_{\mathbb{P}}\left[\iint \mathbf{1}_B(t + s)h(\theta_t, s)T^*(\theta_t, 0, ds)\eta(dt)\right],$$

where we have again used (3.1), this time for the transport $T^*$. By the refined Campbell theorem (2.4),

$$I = \mathbb{E}_{\mathbb{P}_\eta}\left[\iint \mathbf{1}_B(t+s) h(\theta_0, s) T^*(\theta_0, 0, ds)\, dt\right] = \mathbb{E}_{\mathbb{P}_\eta}\left[\int h(\theta_0, s) T^*(\theta_0, 0, ds)\right],$$

where we have used Fubini's theorem and $\lambda(B) = 1$ for the final equation. $\square$

REMARK 3.7. One possible choice of $T$ and $T^*$ in (3.7) is $T(s, \cdot) := \eta$ and $T^*(s, \cdot) := \xi$. Then (3.8) is Neveu's [23] exchange formula

$$(3.9) \qquad \mathbb{E}_{\mathbb{P}_\xi}\left[\int h(\theta_t, -t) \eta(dt)\right] = \mathbb{E}_{\mathbb{P}_\eta}\left[\int h(\theta_0, t) \xi(dt)\right].$$

In case $\xi = \eta$ this is the Mecke equation (2.7).

REMARK 3.8. Let $B, B' \in \mathcal{G}$ have finite and equal Haar measure. Using the definition (2.3) of Palm measures and the invariance of $\xi$ and $\eta$, we can rewrite the exchange formula (3.9) as

$$(3.10) \quad \begin{aligned} &\mathbb{E}_{\mathbb{P}}\left[\iint \mathbf{1}_B(s) h(\theta_t, s-t) \eta(dt) \xi(ds)\right] \\ &= \mathbb{E}_{\mathbb{P}}\left[\iint \mathbf{1}_{B'}(s) h(\theta_s, t-s) \xi(dt) \eta(ds)\right]. \end{aligned}$$

The function $\kappa(\omega, s, t) := h(\theta_s \omega, t - s)$ is invariant in the sense of (3.2). Equation (3.10) implies for all invariant $\kappa$ that

$$(3.11) \quad \begin{aligned} &\mathbb{E}_{\mathbb{P}}\left[\iint \mathbf{1}_B(t) \kappa(s, t) \eta(ds) \xi(dt)\right] \\ &= \mathbb{E}_{\mathbb{P}}\left[\iint \mathbf{1}_{B'}(s) \kappa(s, t) \eta(ds) \xi(dt)\right]. \end{aligned}$$

In case $\xi = \eta$ this gives a version of the mass-transport principle (see [2]) for random measures on Abelian groups. It will be shown in [16] that Neveu's exchange formula (3.9) can be generalized to jointly stationary random measures on a homogeneous space. In fact, the papers [2] and [1] show that the mass-transport principle can be extended beyond this setting.

We finish this section with another useful consequence of Theorem 3.6.

COROLLARY 3.9. *Under the assumption of Theorem 3.6, we have*

$$(3.12) \qquad \mathbb{E}_{\mathbb{P}_\xi}\left[g \int f(\theta_t) T(0, dt)\right] = \mathbb{E}_{\mathbb{P}_\eta}\left[f \int g(\theta_s) T^*(0, ds)\right]$$

*for all measurable functions $f, g: \Omega \to [0, \infty)$.*

PROOF. Apply (3.8) with $h(\omega, s) := f(\omega) g(\theta_s \omega)$. $\square$

**4. Invariance properties of Palm measures.** In this section we fix a stationary $\sigma$-finite measure $\mathbb{P}$ on $(\Omega, \mathcal{F})$. We shall establish fundamental relationships between invariant balancing weighted transport-kernels and Palm measures. Special cases have been known for a long time; see the references below. We were partly motivated by Theorem 16 in the recent paper [11], which deals with $(\lambda, \eta)$-balancing transport-kernels in case $G = \mathbb{R}^d$ and where $\eta$ is an ergodic point process; see Example 4.8.

THEOREM 4.1. *Consider two invariant random measures $\xi$ and $\eta$ on $G$ and an invariant weighted transport-kernel $T$. Then $T$ is $\mathbb{P}$-a.e. $(\xi, \eta)$-balancing iff*

$$\text{(4.1)} \qquad \mathbb{E}_{\mathbb{P}_\xi}\left[\int f(\theta_t) T(0, dt)\right] = \mathbb{E}_{\mathbb{P}_\eta}[f]$$

*holds for all measurable $f : \Omega \to [0, \infty)$.*

For simplicity, we will refer to (4.1) in case $\xi = \eta$ as *invariance* of $\mathbb{P}_\xi$ under $T$.

For the proof of Theorem 4.1, we need the following lemma.

LEMMA 4.2. *Assume that $T$ is a $\mathbb{P}$-a.e. $(\xi, \eta)$-balancing invariant weighted transport-kernel. Then there is an invariant transport-kernel $T^*$ on $G$ such that (3.7) holds for $\mathbb{P}$-a.e. $\omega \in \Omega$.*

PROOF. Similarly as in Example 3.3, we consider the following measure $M$ on $\Omega \times G \times G$:

$$M := \iiint \mathbf{1}\{(\omega, s, t) \in \cdot\} T(\omega, s, dt) \xi(\omega, ds) \mathbb{P}(d\omega).$$

Stationarity of $\mathbb{P}$, (2.2) and (3.1) easily imply that (3.4) holds. Moreover, as (3.3) is assumed to hold for $\mathbb{P}$-a.e. $\omega$, we have

$$\text{(4.2)} \quad M' := \int \mathbf{1}\{(\omega, t) \in \cdot\} M(d(\omega, s, t)) = \iint \mathbf{1}\{(\omega, t) \in \cdot\} \eta(\omega, dt) \mathbb{P}(d\omega).$$

This is a $\sigma$-finite measure on $\Omega \times G$. Similarly as in Example 3.3, we can apply Theorem 3.5 in Kallenberg [13] to obtain an invariant transport-kernel $T^*$ satisfying

$$M = \iint \mathbf{1}\{(\omega, s, t) \in \cdot\} T^*(\omega, t, ds) M'(d(\omega, t)).$$

Recalling the definition of $M$ and the second equation in (4.2), we get for all $A \in \mathcal{F}$ that

$$\mathbb{E}_{\mathbb{P}}\left[\mathbf{1}_A \iint \mathbf{1}\{(s, t) \in \cdot\} T(\theta_0, s, dt) \xi(ds)\right]$$
$$= \mathbb{E}_{\mathbb{P}}\left[\mathbf{1}_A \iint \mathbf{1}\{(s, t) \in \cdot\} T^*(\theta_0, t, ds) \eta(dt)\right].$$

Since $\mathcal{G} \otimes \mathcal{G}$ is countably generated, we obtain the assertion of the lemma. $\square$

PROOF OF THEOREM 4.1. If $T$ is $\mathbb{P}$-a.e. $(\xi, \eta)$-balancing, then (4.1) follows from Lemma 4.2 and Theorem 3.6. Conversely, assume that (4.1) holds. Let $f : \Omega \times G \to [0, \infty)$ be measurable. Using the refined Campbell theorem for $\eta$ and (4.1), we get

$$\mathbb{E}_{\mathbb{P}}\bigg[\int f(\theta_s, s)\eta(ds)\bigg] = \mathbb{E}_{\mathbb{P}_\xi}\bigg[\iint f(\theta_t, s)\, ds\, T(\theta_0, 0, dt)\bigg]$$
$$= \mathbb{E}_{\mathbb{P}_\xi}\bigg[\iint f(\theta_t, s+t) T(\theta_0, 0, dt)\, ds\bigg],$$

where we have used invariance of $\lambda$ and Fubini's theorem for the latter equation. By the refined Campbell theorem for $\xi$ and invariance of $T$, we get that the last term equals

$$\mathbb{E}_{\mathbb{P}}\bigg[\iint f(\theta_{s+t}, s+t) T(\theta_s, 0, dt)\xi(ds)\bigg] = \mathbb{E}_{\mathbb{P}}\bigg[\iint f(\theta_t, t) T(\theta_0, s, dt)\xi(ds)\bigg].$$

Now we combine the latter equations and apply them with $f(\omega, s) := g(\theta_{-s}\omega, s)$, where $g : \Omega \times G \to [0, \infty)$ is measurable. This yields

$$\mathbb{E}_{\mathbb{P}}\bigg[\int g(\theta_0, s)\eta(ds)\bigg] = \mathbb{E}_{\mathbb{P}}\bigg[\iint g(\theta_0, t) T(\theta_0, s, dt)\xi(ds)\bigg].$$

Using this with $g := \mathbf{1}_{A \times B}$, for $A \in \mathcal{F}$ and $B \in \mathcal{G}$, gives

$$\mathbb{E}_{\mathbb{P}}[\mathbf{1}_A \eta(B)] = \mathbb{E}_{\mathbb{P}}\bigg[\mathbf{1}_A \int T(\theta_0, s, B)\xi(ds)\bigg].$$

Since $\mathcal{G}$ is countably generated, this concludes the proof of the theorem. $\square$

EXAMPLE 4.3. Consider the setting of Examples 2.3 and 2.5 and let $\mathbb{P}$ be a $\sigma$-finite stationary measure on $(\Omega, \mathcal{F})$ concentrated on the set $\Omega'$ of all integer-valued $\omega \in \Omega$. Assume that $\mathbb{P}(\tilde{\xi} \in \cdot)$ is $\sigma$-finite. Applying Theorem 4.1 with $T$ given by (3.5) and using (2.5) gives

$$(4.3) \quad \iint \frac{1}{\xi(\omega, \{0\})} \sum_{t : \omega(0, t) > 0} f(\theta_t \omega)\omega(0, t) K(\mu, d\omega) \mathbb{P}_\xi(\tilde{\xi} \in d\mu) = \mathbb{E}_{\mathbb{P}_\eta}[f]$$

for any measurable $f : \Omega \to [0, \infty)$. Specializing to the case of a function $f$ depending only on $\eta(\omega)$ yields Theorem 6.5 in Port and Stone [24]. In the special case $G = \mathbb{R}$ (and under further restrictions on the support of $\mathbb{P}$), (4.3) is Theorem 6.5 in [6].

EXAMPLE 4.4. Let $\xi$ be an invariant (simple) point process on $G$. A *point-allocation for* $\xi$ is an allocation rule $\tau:\Omega \times G \to G$ such that $\tau(s) \in \operatorname{supp}\xi$ whenever $s \in \operatorname{supp}\xi$. Such a point-allocation is called *bijective* if $s \mapsto \tau(s)$ is a bijection on $\operatorname{supp}\xi$ whenever $\xi(G) > 0$. Clearly, this is equivalent to the fact that $\tau$ is $\xi$-preserving; see [28] and [7] for more details. Consider a bijective point-allocation $\tau$ for $\xi$. There is an *inverse* point-allocation $\tau^*$, that is, a bijective point-allocation for $\xi$ satisfying $\tau(\omega, \tau^*(\omega, s)) = \tau^*(\omega, \tau(\omega, s)) = s$ for all $(\omega, s) \in \Omega \times G$ such that $s \in \operatorname{supp}\xi(\omega)$. [Defining $\tau^*(s) := s$ for $s \notin \operatorname{supp}\xi$, it can easily be checked that $\tau^*$ is covariant.] The invariant transport-kernels

$$T(s, \cdot) := \delta_{\tau(s)}, \qquad T^*(s, \cdot) := \delta_{\tau^*(s)}, \qquad s \in G,$$

satisfy (3.7) with $\eta := \xi$. Therefore, we obtain from (3.12) that

$$\mathbb{E}_{\mathbb{P}_\xi}[gf(\theta_\tau)] = \mathbb{E}_{\mathbb{P}_\xi}[fg(\theta_{\tau^*})],$$

where $\theta_\tau : \Omega \to \Omega$ is defined by

(4.4) $$\theta_\tau(\omega) := \theta_{\tau(\omega, 0)}(\omega), \qquad \omega \in \Omega.$$

Taking $g \equiv 1$ yields $\mathbb{E}_{\mathbb{P}_\xi}[f(\theta_\tau)] = \mathbb{E}_{\mathbb{P}_\xi}[f]$, that is, the invariance of $\mathbb{P}_\xi$ under $\theta_\tau$. This is Theorem 3.1 in [7] (cf. also Theorem 9.4.1 in [28]). In fact, this result can also be derived from a more general result in [22].

The results in the previous example can be generalized to invariant random measures $\xi$ and $\eta$. To do so, we consider an allocation rule $\tau$ which is $\mathbb{P}$-a.e. $(\xi, \eta)$-balancing (see Example 3.4) and define the transport-kernel $T$ by $T(s, \cdot) := \delta_{\tau(s)}$. Let $T^*$ be an invariant transport-kernel satisfying (3.7) for $\mathbb{P}$-a.e. $\omega \in \Omega$. Then (3.12) says that

$$\mathbb{E}_{\mathbb{P}_\xi}[gf(\theta_\tau)] = \mathbb{E}_{\mathbb{P}_\eta}\left[f \int g(\theta_s) T^*(0, ds)\right].$$

In particular, we obtain that $\mathbb{P}_\xi(\theta_\tau \in \cdot) = \mathbb{P}_\eta$, where $\theta_\tau$ is defined by (4.4). In case $\xi = \eta$ (and in accordance with the terminology introduced after Theorem 4.1) we will refer to this as *invariance* of $\mathbb{P}_\xi$ *under* $\tau$. Theorem 4.1 implies the following proposition.

PROPOSITION 4.5. *Consider two invariant random measures $\xi$ and $\eta$ and let $\tau$ be an allocation rule. Then $\tau$ is $\mathbb{P}$-a.e. $(\xi, \eta)$-balancing iff*

(4.5) $$\mathbb{P}_\xi(\theta_\tau \in A) = \mathbb{P}_\eta(A), \qquad A \in \mathcal{F}.$$

REMARK 4.6. The invariance of $\mathbb{P}_\xi$ under $\xi$-preserving allocation rules is (essentially) a consequence of Satz 4.3 in [22]. The special case $\xi = \lambda$ was treated in [5].

The *invariant $\sigma$-field* $\mathcal{I} \subset \mathcal{F}$ is the class of all sets $A \in \mathcal{F}$ satisfying $\theta_s A = A$ for all $s \in G$. Let $\xi$ be an invariant random measure with finite intensity and define

$$\hat{\xi} := \mathbb{E}_{\mathbb{P}}[\xi(B)|\mathcal{I}], \tag{4.6}$$

where $B \in \mathcal{G}$ has $\lambda(B) = 1$ and the conditional expectation is defined as for probability measures. (Stationarity implies that this definition is $\mathbb{P}$-a.e. independent of the choice of $B$.) If $\mathbb{P}$ is a probability measure and $G = \mathbb{R}^d$, then $\hat{\xi}$ is called the *sample intensity* of $\xi$; the see [20] and [12]. Assuming that $\mathbb{P}(\hat{\xi} = 0) = 0$, we define the *modified Palm measure* $\mathbb{P}_\xi^*$ of $\xi$ (see [14, 20, 28]) by

$$\mathbb{P}_\xi^*(A) := \mathbb{E}_{\mathbb{P}_\xi}[\hat{\xi}^{-1}\mathbf{1}_A], \qquad A \in \mathcal{F}. \tag{4.7}$$

By this definition and $\hat{\xi} \circ \theta_s = \hat{\xi}$, $s \in G$, we have

$$\mathbb{P}_\xi^*(A) = \mathbb{E}_{\mathbb{P}}\left[\hat{\xi}^{-1} \int \mathbf{1}_A(\theta_s)\mathbf{1}_B(s)\xi(ds)\right] = \mathbb{P}_{\xi'}(A), \qquad A \in \mathcal{F}, \tag{4.8}$$

where the invariant random measure $\xi'$ is defined by $\xi' := \hat{\xi}^{-1}\xi$ if $0 < \hat{\xi} < \infty$ and is the null measure, otherwise. Using (4.8), we obtain the following version of Theorem 4.1.

COROLLARY 4.7. *Consider two invariant random measures $\xi$ and $\eta$ with finite intensities such that $\mathbb{P}(\hat{\xi} = 0) = \mathbb{P}(\hat{\eta} = 0) = 0$ and let $T$ be an invariant weighted transport-kernel. Define $\xi' := \hat{\xi}^{-1}\xi$ and $\eta' := \hat{\eta}^{-1}\eta$. Then $T$ is $\mathbb{P}$-a.e. $(\xi', \eta')$-balancing, iff*

$$\mathbb{E}_{\mathbb{P}_\xi^*}\left[\int f(\theta_t) T(0, dt)\right] = \mathbb{E}_{\mathbb{P}_\eta^*}[f] \tag{4.9}$$

*holds for all measurable $f : \Omega \to [0, \infty)$.*

EXAMPLE 4.8. Let $\eta$ be an invariant random measure with finite intensity and such that $\mathbb{P}(\hat{\eta} = 0) = 0$. Consider the invariant random measure $\xi := \hat{\eta}\lambda$ and let $T$ be an invariant weighted transport-kernel. Then $T$ is $\mathbb{P}$-a.e. $(\xi, \eta)$-balancing iff $T$ is $\mathbb{P}$-a.e. $(\lambda, \eta')$-balancing, where $\eta' := \hat{\eta}^{-1}\eta$. By Corollary 4.7, this is equivalent to

$$\mathbb{E}_{\mathbb{P}}\left[\int \mathbf{1}_A(\theta_t)T(0, dt)\right] = \mathbb{P}_\eta^*(A), \qquad A \in \mathcal{F}. \tag{4.10}$$

In case $G = \mathbb{R}^d$, $\eta$ is a point process, $T$ is Markovian, and $\mathbb{P}$ is an ergodic probability measure; this boils down to Theorem 16 in [11].

EXAMPLE 4.9. Consider Example 4.8 in case $\eta$ is a point process and the weighted transport-kernel is generated by an allocation rule $\tau$ satisfying $\lambda(\{s \in G : \tau(s) \notin \operatorname{supp} \eta\}) = 0$. Clearly, $\tau$ is $\mathbb{P}$-a.e. $(\lambda, \eta')$-balancing iff

$$(4.11) \qquad \lambda(\{s \in G : \tau(s) = t\}) = \hat{\eta}^{-1}, \qquad t \in \operatorname{supp} \eta,$$

holds $\mathbb{P}$-a.e. By Corollary 4.7, this is then equivalent to

$$(4.12) \qquad \mathbb{P}(\theta_\tau \in A) = \mathbb{P}^*_\eta(A), \qquad A \in \mathcal{F}.$$

The special case $G = \mathbb{R}^d$ is Theorem 9.1 in [14], a slight generalization of Theorem 13 in [11]. It is quite remarkable that allocation rules satisfying (4.11) do exist in case $G = \mathbb{R}^d$ (and in case $\mathbb{P}$ is a probability measure); see Theorem 1 in [11] (and Theorem 10.1 in [14] for the nonergodic case). We also refer to the discussion in the introduction and Remark 5.2. Theorem 20 in [11] shows that the situation is different for discrete groups.

REMARK 4.10. Relations (4.5) and (4.12) are examples of *group-coupling* (see [26]); the term "group-coupling" is from [12]. Actually, the relation (4.1) can also be seen as group-coupling by extending the underlying space $(\Omega, \mathcal{F}, \mathbb{P}_\xi)$ to support a random element $\gamma$ in $G$ such that the conditional distribution of $\gamma$ given $\mathcal{F}$ is $T(0, \cdot)$. Then (4.1) can be rewritten as (4.5) with $\theta_\tau$ replaced by $\theta_\gamma$. A similar remark applies to (4.10).

**5. Existence of balancing invariant transport-kernels.** Again we fix a stationary $\sigma$-finite measure $\mathbb{P}$ on $(\Omega, \mathcal{F})$. Our aim is to establish a necessary and sufficient condition for the existence of balancing invariant transport-kernels.

THEOREM 5.1. *Let $\xi$ and $\eta$ be invariant random measures with positive and finite intensities. Then there exists a $\mathbb{P}$-a.e. $(\xi, \eta)$-balancing invariant transport-kernel iff*

$$(5.1) \qquad \mathbb{E}_\mathbb{P}[\xi(B)|\mathcal{I}] = \mathbb{E}_\mathbb{P}[\eta(B)|\mathcal{I}] \qquad \mathbb{P}\text{-}a.e.$$

*for some $B \in \mathcal{G}$ satisfying $0 < \lambda(B) < \infty$.*

PROOF. Let $B \in \mathcal{G}$ satisfy $0 < \lambda(B) < \infty$. For any $A \in \mathcal{I}$, we have from (2.4) that

$$(5.2) \qquad \lambda(B)\mathbb{P}_\xi(A) = \mathbb{E}_\mathbb{P}[\mathbf{1}_A \xi(B)], \qquad \lambda(B)\mathbb{P}_\eta(A) = \mathbb{E}_\mathbb{P}[\mathbf{1}_A \eta(B)].$$

If $T$ is a $\mathbb{P}$-a.e. $(\xi, \eta)$-balancing invariant transport-kernel, then Theorem 4.1 implies the equality $\mathbb{P}_\eta(A) = \mathbb{P}_\xi(A)$ and, thus, $\mathbb{E}_\mathbb{P}[\mathbf{1}_A \eta(B)] = \mathbb{E}_\mathbb{P}[\mathbf{1}_A \xi(B)]$. This entails (5.1).

Let us now assume that (5.1) holds for some $B \in \mathcal{G}$ satisfying $0 < \lambda(B) < \infty$. Since $\mathbb{E}[\xi(\cdot)]$ and $\mathbb{E}[\eta(\cdot)]$ are multiples of $\lambda$, $\xi$ and $\eta$ have the same intensities. We assume without loss of generality that these intensities are equal to 1. From (5.2) and conditioning we obtain that $\mathbb{P}_\xi = \mathbb{P}_\eta$ on $\mathcal{I}$. The group-coupling result in Thorisson [26] (see also Theorem 10.28 in Kallenberg [12]) implies the existence of random elements $\delta$ and $\delta'$ in $\Omega$ and $\rho$ in $G$, all defined on some probability space $(\tilde{\Omega}, \tilde{\mathcal{F}}, \tilde{\mathbb{P}})$, such that $\delta$ has distribution $\mathbb{P}_\xi$, $\delta'$ has distribution $\mathbb{P}_\eta$, and $\delta'(\tilde{\omega}) = \theta_{\rho(\tilde{\omega})} \delta(\tilde{\omega})$ for $\tilde{\mathbb{P}}$-a.e. $\tilde{\omega} \in \tilde{\Omega}$. Let $\tilde{T}(\omega, \cdot)$, $\omega \in \Omega$, be a regular version of the conditional distribution $\tilde{\mathbb{P}}(\rho \in \cdot | \delta = \omega)$. Then we have for any $A \in \mathcal{F}$ that

$$\mathbb{P}_\eta(A) = \tilde{\mathbb{P}}(\delta' \in A) = \tilde{\mathbb{P}}(\theta_\rho \delta \in A) = \mathbb{E}_{\tilde{\mathbb{P}}} \left[ \int \mathbf{1}_A(\theta_s \delta) \tilde{T}(\delta, ds) \right]$$
(5.3)
$$= \mathbb{E}_{\mathbb{P}_\xi} \left[ \int \mathbf{1}_A(\theta_s) \tilde{T}(\theta_0, ds) \right].$$

We now define an invariant transport-kernel $T$ by $T(\omega, s, B) := \tilde{T}(\theta_s \omega, B - s)$. Then (5.3) implies (4.1), and Theorem 4.1 yields that $T$ is $\mathbb{P}$-a.e. $(\xi, \eta)$-balancing. $\square$

REMARK 5.2. The above proof does not provide a method for actually constructing balancing invariant transport-kernels. As mentioned in the Introduction, explicit constructions of allocation rules, in case $\xi = \lambda$ and $\eta$ is an ergodic point process of intensity 1, have been presented by Liggett [19] and Holroyd and Peres [11]. Note that (4.5) means in this case that $\mathbb{P}(\theta_\tau \in \cdot) = \mathbb{P}_\eta$. The construction of balancing invariant transport-kernels and allocation rules in other cases is an interesting topic for further research.

REMARK 5.3. Let $\xi$ and $\eta$ be invariant random measures with finite intensities such that $\mathbb{P}(\hat{\xi} = 0) = \mathbb{P}(\hat{\eta} = 0) = 0$. Then the invariant random measures $\xi' := \hat{\xi}^{-1} \xi$ and $\eta' := \hat{\eta}^{-1} \eta$ satisfy (5.2); see also Corollary 4.7.

**6. Mass-stationarity.** We consider an invariant random measure $\xi$ on $G$ together with a $\sigma$-finite measure $\mathbb{Q}$ on $(\Omega, \mathcal{F})$. Our aim is to establish a condition that is necessary and sufficient for $\mathbb{Q}$ to be the Palm measure of $\xi$ with respect to some stationary $\sigma$-finite measure on $(\Omega, \mathcal{F})$.

Let $C \in \mathcal{G}$ be relatively compact and define an invariant transport-kernel $T_C$ by

(6.1) $\quad T_C(t, B) := \xi(C + t)^{-1} \xi(B \cap (C + t)), \qquad t \in G, \ B \in \mathcal{G},$

if $\xi(C + t) > 0$, and by letting $T_C(t, \cdot)$ equal some fixed probability measure otherwise. In the former case $T_C(t, \cdot)$ is just governing a $G$-valued stochastic experiment that picks a point uniformly in the mass of $\xi$ in $C + t$. If $0 <$ $\lambda(C) < \infty$, we also define the uniform distribution $\lambda_C$ on $G$ by $\lambda_C(B) := \lambda(B \cap C)/\lambda(C)$. The interior (resp. boundary) of a set $C \subset G$ is denoted by $\mathrm{int}\, G$ (resp. $\partial C$).

DEFINITION 6.1. The $\sigma$-finite measure $\mathbb{Q}$ on $(\Omega, \mathcal{F})$ is called *mass-stationary for $\xi$* if $\mathbb{Q}(\xi(G) = 0) = 0$ and

$$(6.2) \quad \mathbb{E}_{\mathbb{Q}}\left[\iint \mathbf{1}_A(\theta_s, s+r)T_C(-r, ds)\lambda_C(dr)\right] = \mathbb{Q} \otimes \lambda_C(A), \qquad A \in \mathcal{F} \otimes \mathcal{G},$$

holds for all relatively compact sets $C \in \mathcal{G}$ with $\lambda(C) > 0$ and $\lambda(\partial C) = 0$.

REMARK 6.2. Assume that $\mathbb{Q}$ is a probability measure. Let $C$ be as assumed in (6.2). Extend the space $(\Omega, \mathcal{F}, \mathbb{Q})$, so as to carry random elements $U, V$ in $G$ such that $\theta_0$ and $U$ are independent, $U$ has distribution $\lambda_C$, and the conditional distribution of $V$ given $(\theta_0, U)$ is uniform in the mass of $\xi$ on $C - U$. (The mappings $\theta_s$, $s \in G$, are extended, so that they still take values in the original space $\Omega$.) Then (6.2) can be written as

$$(6.3) \qquad (\theta_V, U+V) \stackrel{d}{=} (\theta_0, U).$$

In the case of simple point processes on $\mathbb{R}^d$, this is (essentially) the property that was proved in Thorisson ([28], Theorem 9.5.1) to be equivalent to point-stationarity.

THEOREM 6.3. *There exists a $\sigma$-finite stationary measure $\mathbb{P}$ on $(\Omega, \mathcal{F})$ such that $\mathbb{Q} = \mathbb{P}_\xi$ iff $\mathbb{Q}$ is mass-stationary for $\xi$.*

PROOF. Let $C \in \mathcal{G}$ be relatively compact with $\lambda(C) > 0$ and $\lambda(C \setminus \mathrm{int}\, C) = 0$. Then $\lambda(\mathrm{int}\, C) > 0$ and we have for $\lambda$-a.e. $r \in C$ that $r \in \mathrm{int}\, C$. For $r \in \mathrm{int}\, C$ and $t \in G$, we have $t \in \mathrm{int}(C - r + t)$. If, in addition, $t \in \mathrm{supp}\, \xi$, then

$$\xi(C + t - r) \geq \xi(\mathrm{int}(C - r + t)) > 0.$$

By definition (6.1) of $T$ (and using the above fact), we have for all $B, D \in \mathcal{G}$ and $t \in \mathrm{supp}\, \xi$ that

$$\iint \mathbf{1}_B(s)\mathbf{1}_D(s - t + r)T_C(t - r, ds)\lambda_C(dr)$$

$$= \iint \mathbf{1}_B(s)\mathbf{1}_D(s - t + r)\mathbf{1}_{C+t-r}(s)\xi(C + t - r)^{-1}\lambda_C(dr)\xi(ds)$$

$$= \lambda(C)^{-1} \iint \mathbf{1}_B(s)\mathbf{1}_D(r+s)\mathbf{1}_C(r+s)\mathbf{1}_C(r+t)\xi(C - r)^{-1}\, dr\, \xi(ds),$$

where the second equation comes from a change of variables. It follows that

$$\iiint \mathbf{1}_B(s)\mathbf{1}_D(s-t+r)T_C(t-r,ds)\lambda_C(dr)\xi(dt)$$
(6.4)
$$= \lambda(C)^{-1} \iint \mathbf{1}_B(s)\mathbf{1}_D(r+s)\mathbf{1}_C(r+s)\,dr\,\xi(ds)$$
$$= \xi(B)\lambda_C(D).$$

Equation (6.4) implies that the invariant weighted transport-kernel

(6.5) $$T_{C,D}(t,\cdot) := \iint \mathbf{1}\{s \in \cdot\}\mathbf{1}_D(s-t+r)T_C(t-r,ds)\lambda_C(dr)$$

is $(\xi,\eta)$-balancing, where $\eta := \lambda_C(D)\xi$. (Invariance of $T_{C,D}$ is a quick consequence of the same property of $T_C$.) Assume now that $\mathbb{Q} = \mathbb{P}_\xi$ is the Palm measure of $\xi$ with respect to some $\sigma$-finite measure $\mathbb{P}$ on $(\Omega, \mathcal{F})$. Since $\mathbb{P}_\eta = \lambda_C(D)\mathbb{P}_\xi$, we get from Theorem 4.1 (applied with $T = T_{C,D}$) that

$$\mathbb{E}_{\mathbb{P}_\xi}\left[\iint \mathbf{1}_{A'}(\theta_s)\mathbf{1}_D(s+r)T_C(-r,ds)\lambda_C(dr)\right] = \lambda_C(D)\mathbb{P}_\xi(A'), \qquad A' \in \mathcal{F}.$$

This is (6.2) for measurable product sets, implying (6.2) for general $A \in \mathcal{F} \otimes \mathcal{G}$.

Let us now assume, conversely, that $\mathbb{Q}$ is mass-stationary for $\xi$. For simplicity, we can then also assume that $\operatorname{supp}\xi \neq \varnothing$ everywhere on $\Omega$. We will show the Mecke equation (2.7). Let $C \in \mathcal{G}$ be a relatively compact set with $\lambda(C) > 0$ and $\lambda(\partial C) = 0$. Mass-stationarity of $\mathbb{Q}$ implies for any measurable $f: \Omega \to [0,\infty)$ and any $D \in \mathcal{G}$ that

$$\mathbb{E}_{\mathbb{Q}}\left[\iint f(\theta_s)\mathbf{1}_D(s+r)T_C(-r,ds)\lambda_C(dr)\right] = \lambda_C(D)\mathbb{E}_{\mathbb{Q}}[f].$$

By definition (6.1) of $T_C$, this means that

$$\mathbb{E}_{\mathbb{Q}}\left[\iint f(\theta_s)\mathbf{1}_D(s+r)\mathbf{1}_C(s+r)\mathbf{1}_C(r)\xi(C-r)^{-1}\,dr\,\xi(ds)\right]$$
$$= \lambda(D \cap C)\mathbb{E}_{\mathbb{Q}}[f],$$

where we recall the first paragraph of the proof. A change of variables and Fubini's theorem give

$$\int_D \mathbf{1}_C(r)\mathbb{E}_{\mathbb{Q}}\left[\int f(\theta_s)\mathbf{1}_C(r-s)\xi(C-r+s)^{-1}\xi(ds)\right]dr = \mathbb{E}_{\mathbb{Q}}[f]\int_D \mathbf{1}_C(r)\,dr.$$

As $D \in \mathcal{G}$ is arbitrary, this shows that

$$\mathbb{E}_{\mathbb{Q}}\left[\int f(\theta_s)\mathbf{1}_C(r-s)\xi(C-r+s)^{-1}\xi(ds)\right] = \mathbb{E}_{\mathbb{Q}}[f]$$

holds for $\lambda$-a.e. $r \in C$. Applying this with $f$ replaced by $f(\tilde{f} \circ \xi)$, where $\tilde{f}: \mathbf{M} \to [0, \infty)$ is measurable, we obtain for $\lambda$-a.e. $r \in C$ that

$$(6.6) \quad \mathbb{E}_\mathbb{Q}\left[\int f(\theta_s)\tilde{f}(\xi \circ \theta_s)\mathbf{1}_C(r-s)\xi(C-r+s)^{-1}\xi(ds)\right] = \mathbb{E}_\mathbb{Q}[f\tilde{f}(\xi)].$$

By separability of $\mathbf{M}$ (see, e.g., Theorem A2.3 in [12]) and a monotone class argument, we can choose the corresponding null set $C' \in \mathcal{G}$ independently of $\tilde{f}$. Applying (6.6) with $r \in C \setminus C'$ and $\tilde{f}(\mu) := \mu(C-r)$, $\mu \in \mathbf{M}$, gives

$$(6.7) \quad \mathbb{E}_\mathbb{Q}\left[\int f(\theta_s)\mathbf{1}_C(r-s)\xi(ds)\right] = \mathbb{E}_\mathbb{Q}[f\xi(C-r)] \qquad \lambda\text{-a.e. } r \in C.$$

Let $B_n \subset G$, $n \in \mathbb{N}$, be an increasing sequence of compact sets satisfying $\lambda(\partial B_n) = 0$ and $\bigcup_n B_n = G$. (Such a sequence can be constructed with the help of a metric generating the topology on $G$. For any $s \in G$, there is a compact and nonempty ball centred at $s$ whose boundary has $\lambda$-measure 0. Let $B_n^*$, $n \in \mathbb{N}$, be an increasing sequence of compact sets with union $G$. Then $B_n^*$ is contained in the union $\tilde{B}_n$ of finitely many of the above balls. The sequence $\tilde{B}_1 \cup \cdots \cup \tilde{B}_n$, $n \in \mathbb{N}$, has the desired properties.) Fix $n \in \mathbb{N}$ and assume temporarily that

$$(6.8) \quad \mathbb{E}_\mathbb{Q}[f\xi(B_n - B_n)] < \infty,$$

where $B_n - B_n := \{r - r' : r, r' \in B_n\}$. Since $(r, r') \mapsto r - r'$ is continuous, $B_n - B_n$ is again compact. Then we have for all measurable $C' \subset B_n$ and $r \in B_n$ that

$$\mathbb{E}_\mathbb{Q}[f\xi(C'-r)] \leq \mathbb{E}_\mathbb{Q}[f\xi(B_n - B_n)] < \infty.$$

Assume now that $C \subset B_n$ is satisfying the assumptions made in (6.6) and let $C_0 := B_n \setminus C$. Applying (6.7) to $B_n$ yields for $\lambda$-a.e. $r \in B_n$

$$(6.9) \quad \begin{aligned} &\mathbb{E}_\mathbb{Q}\left[\int f(\theta_s)\mathbf{1}_{C_0}(r-s)\xi(ds)\right] + \mathbb{E}_\mathbb{Q}\left[\int f(\theta_s)\mathbf{1}_C(r-s)\xi(ds)\right] \\ &= \mathbb{E}_\mathbb{Q}[f\xi(C_0-r)] + \mathbb{E}_\mathbb{Q}[f\xi(C-r)]. \end{aligned}$$

Since $\partial C_0 \subset \partial B_n \cup \partial(G \setminus C) = \partial B_n \cup \partial C$, we have $\lambda(\partial C_0) = 0$. Hence, we can apply (6.7) to $C_0$ to obtain that the respective first summands in (6.9) coincide $\lambda$-a.e. $r \in C_0$. Therefore, the respective second summands coincide $\lambda$-a.e. $r \in C_0$. Combining this with (6.7) gives

$$(6.10) \quad \mathbb{E}_\mathbb{Q}\left[\int f(\theta_s)\mathbf{1}_C(r-s)\xi(ds)\right] = \mathbb{E}_\mathbb{Q}[f\xi(C-r)], \qquad \lambda\text{-a.e. } r \in B_n.$$

Integrating (6.10) over a measurable set $D \subset B_n$, using (on both sides) Fubini's theorem and a change of variables gives

$$\mathbb{E}_\mathbb{Q}\left[\iint f(\theta_s)\mathbf{1}_D(r+s)\mathbf{1}_C(r)\xi(ds)\,dr\right] = \mathbb{E}_\mathbb{Q}\left[f\iint \mathbf{1}_D(r-s)\mathbf{1}_C(r)\xi(ds)\,dr\right].$$

As both sides are finite measures in $C$ (the right-hand side is bounded by $\mathbb{E}_\mathbb{Q}[f\xi(B_n - B_n)]$) and the class $\mathcal{G}' := \{C \in \mathcal{G} : C \subset B_n, \lambda(\partial C) = 0\}$ is stable under intersections and generates $\mathcal{G} \cap B_n$, we obtain this equation even for all measurable $C \subset B_n$. [To check that $\sigma(\mathcal{G}') = \mathcal{G} \cap B_n$, it is sufficient to show for any nonempty open $U \subset G$ that there is a nonempty open $U' \subset U$ such that $U' \cap B_n \in \mathcal{G}'$. This can be achieved with an open ball $U'$ having $\lambda(\partial U') = 0$.] Reversing the above steps, we obtain (6.10) for all measurable $C \subset B_n$. Since $\mathcal{G}$ is countably generated, we can choose the corresponding null-sets independently of $C$. This means that there is a measurable set $B'_n \subset B_n$ such that $\lambda(B_n \setminus B'_n) = 0$ and

$$(6.11) \quad \mathbb{E}_\mathbb{Q}\left[\int f(\theta_s)\mathbf{1}_C(r-s)\xi(ds)\right] = \mathbb{E}_\mathbb{Q}[f\xi(C-r)], \qquad r \in B'_n, \ C \in \mathcal{G} \cap B_n.$$

Still keeping $n \in \mathbb{N}$ fixed in (6.11), we now lift the assumption (6.8) on $f : \Omega \to [0,\infty)$. If $\mathbb{E}_\mathbb{Q}[f] < \infty$, we can apply (6.11) with $f$ replaced by $f\mathbf{1}\{\xi(B_n - B_n) \le m\}$ and then let $m \to \infty$. For general $f$, we decompose $\Omega$ into measurable sets $D_m \uparrow \Omega$ with $\mathbb{Q}(D_m) < \infty$, apply the previous result to $\mathbf{1}_{D_m} \min\{f, k\}$, and let $m, k \to \infty$. Then (6.11) still holds for all $r \in B''_n \in \mathcal{G}$, where $B''_n \subset B_n$ such that $\lambda(B_n \setminus B''_n) = 0$. For notational simplicity, we assume $B''_n = B'_n$.

In the final step of the proof we would like to take the limit in (6.11) as $n \to \infty$. First we can assume without loss of generality that $\lambda(B_1) > 0$. Let $r_0 \in B_1 \setminus B^*$, where $B^*$ is the $\lambda$-null set $\bigcup_n B_n \setminus B'_n$. Then $r_0 \in B'_n$ for all $n \ge 1$. Take an arbitrary $\tilde{C} \in \mathcal{G}$. Applying (6.11) to $C := \tilde{C} \cap B_n$ and letting $n \to \infty$ yields

$$(6.12) \quad \mathbb{E}_\mathbb{Q}\left[\int f(\theta_s)\mathbf{1}_{C'}(-s)\xi(ds)\right] = \mathbb{E}_\mathbb{Q}\left[\int f\mathbf{1}_{C'}(s)\xi(ds)\right],$$

for $C' = \tilde{C} - r_0$ and hence for any $C' \in \mathcal{G}$. The measure $\mathbb{E}_\mathbb{Q}[\int \mathbf{1}\{(\theta_0, s) \in \cdot\}\xi(ds)]$ is finite on measurable product sets of the form $\{\mu \in D : \mu(B) \le k\} \times B$, where $\mathbb{Q}(D) < \infty$, $B$ is compact, and $k \in \mathbb{N}$. Since $\Omega \times G$ is the monotone union of countably many of such sets, it is now straightforward to proceed from (6.12) to the full Mecke equation (2.7). $\square$

REMARK 6.4. The inversion formula (2.6) implies that the measure $\mathbb{Q}$ in Theorem 6.3 determines $\mathbb{P}$.

Let $C$ be as in (6.2) and assume that $0 \in \mathrm{int}\, C$. Then $\xi(C + t) > 0$ for all $t \in \mathrm{supp}\, \xi$ and one might think (at least at first glance) that a Palm measure of $\xi$ is invariant under $T_C$. The following simple example (other examples can be based on the Poisson process) shows that this is wrong.

EXAMPLE 6.5. Consider the group $G = \{0, 1, 2\}$ with addition modulo 3. Let $\xi_0, \xi_1, \xi_2$ be independent Bernoulli $(1/2)$ random variables. The distribution $\mathbb{P}$ of the point process $\xi_0 \delta_0 + \xi_1 \delta_1 + \xi_2 \delta_2$ is stationary. Let $\mathbb{Q}$ be the Palm probability measure $\mathbb{P}^0_\xi$, defined in the setting of Example 2.1. Since $\lambda$ is (a multiple of) the counting measure, we can take $B := \{0\}$ in (2.3) to see that $\xi\{1\}$ and $\xi\{2\}$ are independent Bernoulli $(1/2)$ under $\mathbb{Q}$. [Of course, we have $\mathbb{Q}(\xi\{0\} = 1) = 1$.] Consider the set $C := \{0, 1\}$ and the event $A := \{\xi\{1\} = 1\}$, where we recall that $\xi$ is the identity on $\Omega = \mathbf{M}$. Then we obtain from a trivial calculation that

$$\mathbb{E}_\mathbb{Q}\left[\int \mathbf{1}_A(\theta_s) T_C(0, ds)\right] = \mathbb{E}_\mathbb{Q}\left[\int \mathbf{1}\{\xi\{1+s\} = 1\} T_C(0, ds)\right] = \frac{3}{8}.$$

Since $\mathbb{Q}(A) = 1/2$, $\mathbb{Q}$ is not invariant under $T_C$.

**7. Discussion of mass-stationarity.** As in Section 6, consider an invariant random measure $\xi$ on $G$ together with a $\sigma$-finite measure $\mathbb{Q}$ on $(\Omega, \mathcal{F})$. Assume that $\mathbb{Q}(\xi(G) = 0) = 0$. In the point process case, [7] and [8] answered a question, raised in [28] and [4], positively by proving that $\mathbb{Q}$ is a Palm measure of $\xi$ if and only if $\mathbb{Q}$ is invariant under all $\sigma(\xi) \otimes \mathcal{G}$-measurable bijective point-allocations $\tau$ for $\xi$, that is, with $\theta_\tau$ defined by (4.4):

(7.1) $$\mathbb{Q}(\theta_\tau \in A) = \mathbb{Q}(A), \qquad A \in \mathcal{F}.$$

At this stage one might be tempted to guess that (7.1) characterizes Palm measures also for general $\xi$. The following example shows that, in general, randomization is needed to define mass-stationarity. We will construct a probability measure $\mathbb{Q}$ and an invariant random measure $\xi$ satisfying (7.1) for all $\xi$-preserving allocation rules $\tau$; see Example 3.4. Still $\mathbb{Q}$ will be no Palm measure of $\xi$. The construction applies to any Abelian group as considered in this paper.

EXAMPLE 7.1. Assume that $(\Omega, \mathcal{F}) = (\mathbf{M} \times \mathbf{M}, \mathcal{M} \otimes \mathcal{M})$ and (abusing notation) $\theta_s \omega = (\theta_s \mu, \theta_s \nu)$ for $\omega = (\mu, \nu) \in \mathbf{M} \times \mathbf{M}$. Let $\Pi$ denote the distribution of a stationary Poisson process with intensity 1, considered as a probability measure on $(\mathbf{M}, \mathcal{M})$. It is well known that the associated Palm probability measure (defined in the framework of Example 2.1) is given by $\Pi^0 = \int \mathbf{1}\{\mu + \delta_0 \in \cdot\} \Pi(d\mu)$. Let $\mathbb{Q} := \Pi^0 \otimes \Pi$ and $c > 0$ be an irrational number. Define an invariant random measure $\xi$ on $G$ by $\xi := \xi_1 + c\xi_2$, where $\xi_1$ and $\xi_2$ are the projections of $\Omega$ onto the first and second component, respectively. Let $\tau$ be a $\xi$-preserving allocation rule, that is,

$$\int \mathbf{1}\{\tau(s) = t\} \xi_1(ds) + c \int \mathbf{1}\{\tau(s) = t\} \xi_2(ds) = \xi\{t\}, \qquad t \in G.$$

Since $c$ is irrational, this can only hold if $\tau$ is $\xi_1$-preserving. As it can be straightforwardly checked that $\mathbb{Q}$ is the Palm measure of $\xi_1$, Theorem 4.1 implies (7.1).

We now show that $\mathbb{Q}$ is not mass-stationary for $\xi$. Therefore and by Theorem 6.3, it cannot be a Palm measure of $\xi$. Consider a set $C \in \mathcal{G}$ as in Definition 6.1 and place it at random around the origin. This random set will contain a $\xi_2$-point with positive probability; this point will in turn be chosen with positive probability as a new origin. Thus, if $\mathbb{Q}$ was mass-stationary for $\xi$, it should also have mass $c$ at 0 with positive probability. But this is not the case.

Instead of working with general invariant (weighted) transport-kernels, we define mass-stationarity by (6.2). This property has the advantage of having the direct probabilistic interpretation (6.3) when $\mathbb{Q}$ is a probability measure. In order to see how it is related to invariance under weighted transport-kernels, let $C \in \mathcal{G}$ be as in Definition 6.1 and $D \in \mathcal{G}$ with $\lambda(C \cap D) > 0$. Define the invariant weighted transport-kernel $T'_{C,D} := \lambda_C(D)^{-1} T_{C,D}$, where $T_{C,D}$ is given by (6.4). As noted at (6.4), we have that $T'_{C,D}$ is $\xi$-preserving. Also $T'_{C,D}$ is bounded but in general not Markovian. Mass-stationarity of $\mathbb{Q}$ is equivalent to assuming invariance of $\mathbb{Q}$ under all these transport-kernels. Now, Theorem 6.1 and Theorem 4.1 yield the following result.

THEOREM 7.2. *The measure $\mathbb{Q}$ is mass-stationary for $\xi$ iff it is invariant under bounded $\xi$-preserving invariant weighted transport-kernels $T$.*

Further results on mass-stationary random measures will be provided in the papers [15, 17] and [18].

We finish this section with some open problems related to mass-stationarity. A kernel $T$ from $\Omega \times G$ to $G$ is called $\xi$-*measurable* if $T(\cdot, \cdot, B)$ is $\sigma(\xi) \otimes \mathcal{G}$-measurable for all $B \in \mathcal{G}$.

PROBLEM 7.3. Assume that $\mathbb{Q}$ is invariant under $\xi$-preserving and $\xi$-measurable invariant transport-kernels. Is $\mathbb{Q}$ mass-stationary for $\xi$?

The condition in Problem 7.3 implies that of the following problem.

PROBLEM 7.4. Assume that

$$(7.2) \quad \mathbb{E}_\mathbb{Q}\left[\iint \mathbf{1}_A(\theta_s) T_C(-r, ds) \lambda_C(dr)\right] = \mathbb{Q}(A), \qquad A \in \mathcal{F},$$

holds for all $C$ as in Definition 6.1. Is $\mathbb{Q}$ mass-stationary for $\xi$?

The counterexample in Example 7.1 arises because mass-atoms of relatively prime size cannot be mapped into each other in a measure-preserving way. But what about diffuse random measures?

PROBLEM 7.5. Assume that $\xi$ is diffuse and that (7.1) holds for all $\xi$-preserving [and $\sigma(\xi) \otimes \mathcal{G}$-measurable] allocation rules $\tau$. Is $\mathbb{Q}$ mass-stationary for $\xi$?

If the answer to Problem 7.5 is negative, we might (as is the key idea in [27] and [28]) attempt to introduce a *stationary independent background*:

PROBLEM 7.6. Let $\theta_t$, $t \in G$, and $\xi$ be defined on $(\Omega, \mathcal{F}, \mathbb{Q})$. Introduce a *stationary independent background* as follows: let $\theta'_t$, $t \in G$, be another flow defined on a space $(\Omega', \mathcal{F}', \mathbb{Q}')$, where $\mathbb{Q}'$ is stationary under the flow, and consider the joint flow on $(\Omega, \mathcal{F}, \mathbb{Q}) \otimes (\Omega', \mathcal{F}', \mathbb{Q}')$ with $\xi$ defined in the natural way on this extended space. Assume that $\xi$ is diffuse and that $\mathbb{Q} \otimes \mathbb{Q}'$ is invariant under $\xi$-preserving [and $\sigma(\xi) \otimes \mathcal{G}$-measurable] allocation rules for all such stationary independent backgrounds. Is $\mathbb{Q}$ mass-stationary for $\xi$?

PROBLEM 7.7. Same as Problem 7.6 but now only assume that $\mathbb{Q}$—and not $\mathbb{Q} \otimes \mathbb{Q}'$—is invariant.

REMARK 7.8. The second author wants to use this oportunity to correct a mistake in [27] and [28]. Remark 3.2 in [27] claims that the answer to Problem 7.7 is positive in the case of simple point processes on $\mathbb{R}^d$. This claim is a mistake stemming from the author's forgetting that the argument in Section 4.4 in [27] relies on the joint invariance. Similarly, Lemma 4.1 in [27] needs to be corrected by adding the background in (4.11). The same applies to Remark 9.3.2 and Lemma 9.4.1 in [28].

**Acknowledgments.** We would like to thank Alexander Holroyd, Yuval Peres and an anonymous referee for valuable comments.

## REFERENCES

[1] ALDOUS, D. and LYONS, R. (2007). Processes on unimodular random networks. *Electron. J. Probab.* **12** 1454–1508 (electronic). MR2354165
[2] BENJAMINI, I., LYONS, R., PERES, Y. and SCHRAMM, O. (1999). Group-invariant percolation on graphs. *Geom. Funct. Anal.* **9** 29–66. MR1675890
[3] CHATTERJEE, S., PELED, R., PERES, Y. and ROMIK, R. (2008). Gravitational allocation to Poisson points. *Ann. Math.* To appear.
[4] FERRARI, P. A., LANDIM, C. and THORISSON, H. (2004). Poisson trees, succession lines and coalescing random walks. *Ann. Inst. H. Poincaré Probab. Statist.* **40** 141–152. MR2044812

[5] GEMAN, D. and HOROWITZ, J. (1975). Random shifts which preserve measure. *Proc. Amer. Math. Soc.* **49** 143–150. MR0396907
[6] HARRIS, T. E. (1971). Random measures and motions of point processes. *Z. Wahrsch. Verw. Gebiete* **18** 85–115. MR0292148
[7] HEVELING, M. and LAST, G. (2005). Characterization of Palm measures via bijective point-shifts. *Ann. Probab.* **33** 1698–1715. MR2165576
[8] HEVELING, M. and LAST, G. (2007). Point shift characterization of Palm measures on Abelian groups. *Electron. J. Probab.* **12** 122–137 (electronic). MR2280261
[9] HOLROYD, A. E. and LIGGETT, T. M. (2001). How to find an extra head: Optimal random shifts of Bernoulli and Poisson random fields. *Ann. Probab.* **29** 1405–1425. MR1880225
[10] HOLROYD, A. E. and PERES, Y. (2003). Trees and matchings from point processes. *Electron. Comm. Probab.* **8** 17–27 (electronic). MR1961286
[11] HOLROYD, A. E. and PERES, Y. (2005). Extra heads and invariant allocations. *Ann. Probab.* **33** 31–52. MR2118858
[12] KALLENBERG, O. (2002). *Foundations of Modern Probability*, 2nd ed. Springer, New York. MR1876169
[13] KALLENBERG, O. (2007). Invariant measures and disintegrations with applications to Palm and related kernels. *Probab. Theory Related Fields* **139** 285–310. MR2322698
[14] LAST, G. (2006). Stationary partitions and Palm probabilities. *Adv. in Appl. Probab.* **38** 602–620. MR2256871
[15] LAST, G. (2009). Modern random measures: Palm theory and related models. In *New Perspectives in Stochastic Geometry* (W. Kendall und I. Molchanov, eds.). Clarendon Press, Oxford. To appear.
[16] LAST, G. (2009). Stationary random measures on homogeneous spaces. To appear.
[17] LAST, G. and THORISSON, H. (2009). Characterization of mass-stationary by Bernoulli and Cox transports. To appear.
[18] LAST, G. and THORISSON, H. (2008). Constructions of stationary and mass-stationary random measures. (In preparation.)
[19] LIGGETT, T. M. (2002). Tagged particle distributions or how to choose a head at random. In *In and Out of Equilibrium (Mambucaba, 2000). Progr. Probab.* **51** 133–162. Birkhäuser, Boston. MR1901951
[20] MATTHES, K., KERSTAN, J. and MECKE, J. (1978). *Infinitely Divisible Point Processes*. Wiley, Chichester. MR0517931
[21] MECKE, J. (1967). Stationäre zufällige Masse auf lokalkompakten Abelschen Gruppen. *Z. Wahrsch. Verw. Gebiete* **9** 36–58. MR0228027
[22] MECKE, J. (1975). Invarianzeigenschaften allgemeiner Palmscher Maße. *Math. Nachr.* **65** 335–344. MR0374385
[23] NEVEU, J. (1977). Processus ponctuels. In *École D'Été de Probabilités de Saint-Flour, VI—1976. Lecture Notes in Mathematics* **598** 249–445. Springer, Berlin. MR0474493
[24] PORT, S. C. and STONE, C. J. (1973). Infinite particle systems. *Trans. Amer. Math. Soc.* **178** 307–340. MR0326868
[25] RACHEV, S. T. and RÜSCHENDORF, L. (1998). *Mass Transportation Problems. Vol. I: Theory*. Springer, New York. MR1619170
[26] THORISSON, H. (1996). Transforming random elements and shifting random fields. *Ann. Probab.* **24** 2057–2064. MR1415240
[27] THORISSON, H. (1999). Point-stationarity in *d* dimensions and Palm theory. *Bernoulli* **5** 797–831. MR1715440

[28] THORISSON, H. (2000). *Coupling, Stationarity, and Regeneration.* Springer, New York. MR1741181
[29] TIMÁR, Á. (2004). Tree and grid factors for general point processes. *Electron. Comm. Probab.* **9** 53–59 (electronic). MR2081459
[30] ZÄHLE, U. (1988). Self-similar random measures. I. Notion, carrying Hausdorff dimension, and hyperbolic distribution. *Probab. Theory Related Fields* **80** 79–100. MR970472

INSTITUT FÜR STOCHASTIK
UNIVERSITÄT KARLSRUHE
76128 KARLSRUHE
GERMANY
E-MAIL: last@math.uni-karlsruhe.de

SCIENCE INSTITUTE
UNIVERSITY OF ICELAND
DUNHAGA 3
107 REYKJAVIK
ICELAND
E-MAIL: hermann@hi.is